\def \sgn{\mathop{\rm sgn} \nolimits}
\def\S{{\cal S}}
\baselineskip=14pt
\parskip=10pt
\def\Tilde{\char126\relax}
\font\eightrm=cmr8  
\font\eighttt=cmtt8
\magnification=\magstephalf
\parindent=0pt
\overfullrule=0in
 
\bf
\centerline
{Aufgabe VII.47 of P\'olya-Szeg\"o Immediately Implies 
Dave Robbins's Multi-Integral Evaluation }
\rm
\bigskip
\centerline{ {\it Doron ZEILBERGER}\footnote{$^1$}
{\eightrm  \raggedright
Department of Mathematics, Temple University,
Philadelphia, PA 19122, USA. 
{\eighttt zeilberg@math.temple.edu} \hfill \break
{\eighttt http://www.math.temple.edu/\Tilde zeilberg   .}
May 28, 1998.  Exclusive to the author's website and the
xxx archives.
} 
}
 
Exercise VII.47 of [PS] (brought to my attention by Richard Stanley),
$$
 \sum_{\pi \in \S_k} \sgn (\pi) \cdot
\pi \left [ 
{
{x_1 x_2^2 \dots x_k^k}
\over
{(1- x_k)(1- x_k x_{k-1} ) \cdots (1 - x_k x_{k-1} \dots x_1 ) }
}
\right ]  =
{
{x_1 \cdots x_k \prod_{1 \leq i < j \leq k} ( x_j - x_i )}
\over
{ \prod_{i=1}^k (1- x_i ) \prod_{1 \leq i < j \leq k} (1- x_i x_j ) }
} \quad ,
$$
(that is easily proved by induction on $k$ and Lagrange Interpolation),
immediately implies the main result of [R], upon
setting $x_i:=q^{a_i}$, multiplying by $(1-q)^k$, and 
letting $q \rightarrow 1$.
 
{\bf References}
 
[PS] George P\'olya and Gabor Szeg\"o, 
{\it Aufgaben und Lehrs\"atze aus der Analysis}, Julius Springer,
Berlin, 1925.
 
[R] David Robbins, {\it An application of Okada's Minor Summation
Formula to the Evaluation of a Multiple Integral},
xxx archives ({\tt http://xxx.lanl.gov}), {\tt math.CO/9805108},
23 May 1998, (3 pages).
 
\bye